\documentclass[a4paper]{amsart}
\usepackage{amsmath,amsthm,amssymb,latexsym,epic,bbm,comment}
\usepackage{graphicx,enumerate,stmaryrd}
\usepackage[all,2cell]{xy}
\xyoption{2cell}

\newtheorem{theorem}{Theorem}
\newtheorem{lemma}[theorem]{Lemma}
\newtheorem{corollary}[theorem]{Corollary}
\newtheorem{proposition}[theorem]{Proposition}

\usepackage[all]{xy}
\usepackage[active]{srcltx}
\usepackage[parfill]{parskip}
\usepackage{enumerate}

\font\sc=rsfs10
\newcommand{\cC}{\sc\mbox{C}\hspace{1.0pt}}

\font\scc=rsfs7
\newcommand{\ccC}{\scc\mbox{C}\hspace{1.0pt}}

\begin{document}
\title[Characterisation and applications of $\Bbbk$-split bimodules]
{Characterisation and applications\\ of $\Bbbk$-split bimodules}

\author{Volodymyr Mazorchuk, Vanessa Miemietz and Xiaoting Zhang}

\begin{abstract}
We describe the structure of bimodules (over finite dimensional algebras)
which have the property that the functor of tensoring with such a bimodule
sends any module to a projective module. The main result is that all
such bimodules are $\Bbbk$-split in the sense that they factor 
(inside the tensor category of bimodules) over $\Bbbk$-vector spaces.
As one application, we show that any simple $2$-category
has a faithful $2$-representation inside the $2$-category of 
$\Bbbk$-split bimodules. As another application, we classify 
simple transitive $2$-rep\-re\-sen\-ta\-ti\-ons of the $2$-category of 
projective bimodules over the algebra $\Bbbk[x,y]/(x^2,y^2,xy)$.
\end{abstract}

\maketitle

\section{Introduction and description of the results}\label{s1}

The structure and representation theory of $2$-categories is a young and intensively
studied area of modern mathematics which originated in \cite{BFK,Kh,CR,Ro,KL}.
The series of papers \cite{MM1}--\cite{MM6} initiated the study of the structure 
and representation theory of so-called {\em finitary $2$-categories}, which
are natural $2$-analogues of finite dimensional associative algebras.

Classical representation theory of finite dimensional algebras is essentially based on the classification of simple algebras provided by the classical Artin-Wedderburn Theorem.
For a special class of finitary $2$-categories, called {\em strongly regular fiat $2$-categories},
an analogue of the Artin-Wedderburn Theorem was obtained in \cite{MM3}. Fiat $2$-categories
can be viewed as a vague $2$-analogue of cellular algebras: they have a weak involution
and adjunction morphisms which lead to the fact that, in each $2$-representation, the
involution gets interpreted as taking the adjoint functor. This involution plays 
a crucial role in all arguments of \cite{MM3} and therefore it is unlikely that 
one could generalise these arguments to any wider class of $2$-categories.

In the present paper we take a first step towards understanding the structure of 
arbitrary simple finitary $2$-categories. The main idea motivating our study 
comes, fairly unexpectedly, from the results of \cite{KM2,KMMZ}. Following the ideas
in the proof of \cite[Theorem~2]{KMMZ}, which are based on the main result of \cite{KM2},
we observe that every simple finitary $2$-category can be faithfully represented 
using functorial actions in which all involved functors are right exact and
have the property that they send any module to a projective module.

The main technical result of the this article is a characterisation and description
of such functors. In fact, in Theorem~\ref{thm1} we show that bimodules representing such
functors belong to the additive closure of bimodules of the form $P\otimes_{\Bbbk}N$,
where $P$ is a projective left module and $N$ is an arbitrary right module. Put differently,
in the tensor category of bimodules, all such bimodules factor through $\Bbbk$-vector spaces.
A precise formulation of Theorem~\ref{thm1} and its proof can be found in Section~\ref{s3}.

We also give two applications of Theorem~\ref{thm1}. The first one, which can be found in 
Section~\ref{s4}, concerns the faithful representation of simple finitary $2$-categories
as explained above. The second one, presented in Section~\ref{s5}, is of different nature.
It concerns the problem of classifying simple transitive $2$-representations for the
$2$-category of projective bimodules over the finite dimensional algebra 
$A:=\Bbbk[x,y]/(x^2,y^2,xy)$. This kind of problem was studied, in different contexts, for many 
other $2$-categories, see \cite{MM5,MM6,Zi,MaMa,KMMZ,MT,MMMT,MZ,Zh2}. 
The importance of this problem is supported by interesting applications, see e.g. \cite{KM1}.
The case of the algebra $A$ treated in this paper differs
significantly from all previously studied cases. To start with, the $2$-category of
projective functors for $A$ is not fiat, as $A$ is not self-injective, cf. 
\cite[Subsection~7.3]{MM1}. However, simple transitive $2$-representations for some non-fiat $2$-categories have also been classified in \cite{Zh2,MZ}. The main point of 
the algebra $A$ is that this is the smallest algebra which does not have any non-zero
projective-injective modules. Therefore the general approach outlined in \cite{MZ},
which is based on existence of a projective-injective module, is
not applicable either. In Section~\ref{s5} we propose a new approach to this problem
which crucially uses our Theorem~\ref{thm1}. In Section~\ref{snew} we show that
the decategorification of a fiat $2$-category with strongly regular two-sided cells
is a quasi-hereditary algebra with a simple preserving duality.

As the material discussed in Sections~\ref{s3}, \ref{s4} and \ref{s5} is of rather different
nature, we do not provide a general list of notation and preliminaries but rather 
postpone all this to each individual section separately.

For some further examples and structural results on finitary $2$-categories we refer 
the reader to \cite{GM1,GM2,Xa,Zh1}.
\vspace{0.5cm}
 
\textbf{Acknowledgements:} This research was partially supported by
the Swedish Research Council, Knut and Alice Wallenberg Stiftelse and 
G{\"o}ran Gustafsson Stiftelse. We thank Claus Michael Ringel for stimulating discussions.

\section{$\Bbbk$-split bimodules for finite dimensional algebras}\label{s3}

\subsection{Main Theorem}\label{s3.1}

Throughout the paper, we fix an algebraically closed field $\Bbbk$. For finite dimensional
$\Bbbk$-algebras $A$ and $B$, we denote by
\begin{itemize}
\item $A$-mod the category of finite dimensional left $A$-modules;
\item mod-$A$ the category of finite dimensional right $A$-modules;
\item $A$-mod-$B$ the category of finite dimensional $A$-$B$-bimodules.
\end{itemize}

The main result of the paper is the following statement.

\begin{theorem}\label{thm1}
Let $A$ and $B$ be two basic finite dimensional $\Bbbk$-algebras 
and $Q\in A\text{-}\mathrm{mod}\text{-}B$. Then the following conditions are equivalent.
\begin{enumerate}[$($a$)$]
\item\label{thm1.1} The functor $Q\otimes_B{}_-:B\text{-}\mathrm{mod}\to
A\text{-}\mathrm{mod}$ maps any $B$-module to a projective $A$-module.
\item\label{thm1.2} The functor $\mathrm{Hom}_{A\text{-}}(Q,{}_-):A\text{-}\mathrm{mod}\to
B\text{-}\mathrm{mod}$ maps any short exact sequence in $A\text{-}\mathrm{mod}$
to a split short exact sequence in $B\text{-}\mathrm{mod}$.
\item\label{thm1.3} The $A$-$B$-bimodule $Q$ belongs to the additive closure,
in $A\text{-}\mathrm{mod}\text{-}B$, of all
$A$-$B$-bimodules of the form $A\otimes_{\Bbbk}K$, where $K\in\mathrm{mod}\text{-}B$.
\end{enumerate}
\end{theorem}

Bimodules satisfying the conditions in Theorem~\ref{thm1}\eqref{thm1.3} will be called {\em $\Bbbk$-split}.

Note that, by considering the algebra $A\times B$, we can reduce Theorem~\ref{thm1}
to the case $A=B$. So, in the proof which follows, we assume $A=B$.

\subsection{Implication \eqref{thm1.1}$\Rightarrow$\eqref{thm1.2}}\label{s3.2}

Assume that condition~\eqref{thm1.1} is satisfied. Then, in particular,
$Q\otimes_A A\cong {}_AQ$ is projective and hence the functor 
$\mathrm{Hom}_{A\text{-}}(Q,{}_-)$ is exact. Furthermore, for any $M\in A\text{-}\mathrm{mod}$,
the $A$-module $Q\otimes_A M$ is projective and therefore the functor
\begin{displaymath}
\mathrm{Hom}_{A\text{-}}(Q\otimes_A M,{}_-)\cong
\mathrm{Hom}_{A\text{-}}(M,\mathrm{Hom}_{A\text{-}}(Q,{}_-)):A\text{-}\mathrm{mod}\to
\Bbbk\text{-}\mathrm{mod} 
\end{displaymath}
is exact.

For a short exact sequence
\begin{displaymath}
0\to X\to Y\to Z \to 0
\end{displaymath}
in $A\text{-}\mathrm{mod}$, application of the exact functor
$\mathrm{Hom}_{A\text{-}}(Q,{}_-)$ produces the short exact sequence
\begin{equation}\label{eq1}
0\to \mathrm{Hom}_{A\text{-}}(Q,X)\to\mathrm{Hom}_{A\text{-}}(Q,Y)\overset{\alpha}{\to}\mathrm{Hom}_{A\text{-}}(Q,Z)\to 0.
\end{equation}
If \eqref{eq1} splits, then its image under
$\mathrm{Hom}_{A\text{-}}(M,{}_-)$ is, clearly,  split short exact, for any $M\in A\text{-}\mathrm{mod}$.
At the same time, if \eqref{eq1} does not split, then,
for $M=\mathrm{Hom}_{A\text{-}}(Q,Z)$, the identity morphism on $M$ is not in the
image of the map
\begin{displaymath}
\mathrm{Hom}_{A\text{-}}(M,\mathrm{Hom}_{A\text{-}}(Q,Y))\overset{\alpha\circ{}_-}{\longrightarrow}
\mathrm{Hom}_{A\text{-}}(M,M).
\end{displaymath}
Therefore the latter map is not surjective and, consequently, the sequence
\begin{displaymath}
0\to \mathrm{Hom}_{A\text{-}}(M,\mathrm{Hom}_{A\text{-}}(Q,X))\to
\mathrm{Hom}_{A\text{-}}(M,\mathrm{Hom}_{A\text{-}}(Q,Y))\to
\mathrm{Hom}_{A\text{-}}(M,M)\to 0
\end{displaymath}
is not exact. Thus the functor $\mathrm{Hom}_{A\text{-}}(Q\otimes M,{}_-)$
is not exact either, a contradiction. Hence condition~\eqref{thm1.2} is satisfied.

\subsection{Implication \eqref{thm1.2}$\Rightarrow$\eqref{thm1.3}}\label{s3.3}

Assume that condition~\eqref{thm1.2} is satisfied. In particular, the functor
$\mathrm{Hom}_{A\text{-}}(Q,{}_-)$ is exact and thus the left $A$-module ${}_{A}Q$ is projective.
Denote by $R$ the Jacobson radical $\mathrm{Rad}(A)$ of $A$.

Applying the duality $*:=\mathrm{Hom}_{\Bbbk}({}_-,\Bbbk)$ to the short exact sequence
\begin{displaymath}
0\to R\to A\to \mathrm{top}(A)\to 0 
\end{displaymath}
in $A\text{-}\mathrm{mod}\text{-}A$, gives the short exact sequence
\begin{displaymath}
0\to (\mathrm{top}(A))^* \to A^*\to R^*\to 0
\end{displaymath}
in $A\text{-}\mathrm{mod}\text{-}A$. Applying $\mathrm{Hom}_{A\text{-}}(Q,{}_-)$ 
to the latter short exact sequence results in the short exact sequence
\begin{equation}\label{eq2}
0\to\mathrm{Hom}_{A\text{-}}(Q,(\mathrm{top}(A))^*)\to
\mathrm{Hom}_{A\text{-}}(Q,A^*)\to\mathrm{Hom}_{A\text{-}}(Q,R^*)\to 0.
\end{equation}
By condition~\eqref{thm1.2}, this sequence is split in $A\text{-}\mathrm{mod}$.

By adjunction, we have
\begin{equation}\label{eq3}
\mathrm{Hom}_{A\text{-}}(Q,A^*)\cong\mathrm{Hom}_{\Bbbk}(A\otimes_A Q,\Bbbk)\cong Q^*
\end{equation}
and
\begin{equation}\label{eq4}
\mathrm{Hom}_{A\text{-}}(Q,R^*)\cong \mathrm{Hom}_{\Bbbk}(R\otimes_AQ,\Bbbk). 
\end{equation}
Moreover, as any $A$-homomorphism from $Q$ to $(\mathrm{top}(A))^*$ 
vanishes on $RQ$, we have
\begin{displaymath}
\mathrm{Hom}_{A\text{-}}(Q,(\mathrm{top}(A))^*)\cong 
\mathrm{Hom}_{A\text{-}}(Q/RQ,(\mathrm{top}(A))^*) 
\end{displaymath}
and then, by adjunction,
\begin{displaymath}
\mathrm{Hom}_{A\text{-}}(Q/RQ,(\mathrm{top}(A))^*)\cong
\mathrm{Hom}_{\Bbbk}(\mathrm{top}(A)\otimes_A (Q/RQ),\Bbbk).
\end{displaymath}
Finally, we observe that $Q/RQ$ is semi-simple as left $A$-module which
yields an isomorphism $\mathrm{top}(A)\otimes_A(Q/RQ)\cong Q/RQ$.

Plugging the latter into \eqref{eq2}, using \eqref{eq3} and \eqref{eq4}, 
and applying $*$, gives us the short exact sequence
\begin{displaymath}
0\to R\otimes_AQ\to Q\to Q/RQ\to 0
\end{displaymath}
which is split in $\mathrm{mod}\text{-}A$. We denote by
$\beta:Q/RQ \to Q$ the splitting morphism in $\mathrm{mod}\text{-}A$.

Fix a decomposition of the $A$-$A$-bimodule $Q/RQ$ into a direct sum
$X_1\oplus X_2\oplus \dots\oplus X_k$ of indecomposable
$A$-$A$-bimodules. As $Q/RQ$ is semi-simple, as a left $A$-module,
we have that each $X_i$ is isotypic as a left $A$-module and indecomposable
as a right $A$-module. Let $L_i$ denote the (unique) simple subquotient of 
the left $A$-module ${}_AX_i$ and $P_i$ denote an indecomposable
projective cover of $L_i$ in $A\text{-}\mathrm{mod}$. Let also $e_i$ denote
a primitive idempotent corresponding to $P_i$.

Consider the $A$-$A$-bimodule
\begin{displaymath}
\hat{Q}:=\bigoplus_{i=1}^k P_i\otimes_{\Bbbk} X_i 
\end{displaymath}
and note that $\hat{Q}$ belongs to the additive closure of bimodules described
in condition~\eqref{thm1.3}. By adjunction, for any $A$-$A$-bimodule $V$, we have
\begin{equation}\label{eq5}
\mathrm{Hom}_{A\text{-}A}(\hat{Q},V)\cong
\bigoplus_{i=1}^k\mathrm{Hom}_{\text{-}A}(X_i,e_iV).
\end{equation}

The homomorphism $\beta$ induces, by equation~\eqref{eq5}, a homomorphism $\gamma:\hat{Q}\to Q$
of $A$-$A$-bimodules. By construction, the image of this homomorphism covers
$Q/RQ$, that is the top of $Q$, considered as a left $A$-module. Therefore
$\gamma$ is surjective by the Nakayama Lemma. At the same time, we already know that
$Q$ is projective as a left $A$-module and that all $X_i$ are isotypic as left
$A$-modules. Compared to the construction of $\hat{Q}$, this 
implies the equality $\dim(\hat{Q})=\dim(Q)$ and yields that $\gamma$
is, in fact, an isomorphism. Condition~\eqref{thm1.3} follows.

\subsection{Implication \eqref{thm1.3}$\Rightarrow$\eqref{thm1.1}}\label{s3.4}

Assume that condition~\eqref{thm1.3} is satisfied. Note that, for any 
$M\in A\text{-}\mathrm{mod}$, the $A$-module $A\otimes_{\Bbbk} K\otimes_AM$ is just 
a direct sum of $\dim(K\otimes_A M)$ copies of $A$, in particular, it
is projective. By additivity, $Q\otimes_A M$ is also projective, for any $Q$
from the additive closure of the bimodules described in condition~\eqref{thm1.3}.
Therefore condition~\eqref{thm1.1} is satisfied.

\section{Application: simple finitary $2$-categories}\label{s4}

\subsection{Finitary $2$-categories}\label{s4.1}

For generalities on $2$-categories we refer to \cite{McL,Le}.
Following \cite[Subsection~2.2]{MM1}, a {\em finitary $2$-category over $\Bbbk$} is a $2$-category 
$\cC$ such that
\begin{itemize}
\item $\cC$ has finitely many objects;
\item each $\cC(\mathtt{i},\mathtt{j})$ is equivalent to the category of projective modules
over some finite dimensional associative $\Bbbk$-algebra
(i.e. is a {\em finitary $\Bbbk$-linear} category);
\item all compositions are biadditive and $\Bbbk$-bilinear, when applicable;
\item identity $1$-morphisms are indecomposable.
\end{itemize}

In particular, we have a finite set $\mathcal{S}[\cC]$ of isomorphism classes of indecomposable
$1$-morphisms. Furthermore, by \cite[Section~3]{MM2}, $\mathcal{S}[\cC]$ has the natural structure
of a multisemigroup (cf. \cite{KuM}). 
We denote by $\leq_L$ and $\sim_L$ the corresponding 
{\em left order} and {\em left equivalence relation}, by $\leq_R$ and $\sim_R$ the corresponding 
{\em right order} and {\em right equivalence relation} and by $\leq_J$ and $\sim_J$ the corresponding 
{\em two-sided order} and {\em two-sided equivalence relation}, see \cite[Section~3]{MM2}
and \cite[Subsection~4.1]{KuM}. Equivalence classes for $\sim_L$, $\sim_R$ and $\sim_J$ are
called {\em cells} (left, right or two-sided, respectively).

A $2$-category which satisfies all the above conditions except for the last one will be called
{\em weakly finitary}. Clearly, splitting idempotents in the endomorphism algebras of the
identity $1$-morphisms, every weakly finitary $2$-category can be Morita reduced to a
finitary $2$-category (cf. \cite{MM4} for general Morita theory of finitary $2$-categories).

A finitary $2$-category $\cC$ will be called {\em simple} provided that 
\begin{itemize}
\item any non-zero $2$-ideal of $\cC$ contains the identity $2$-morphism for some
non-zero $1$-morphism;
\item there is a unique two-sided cell, which we call $\mathcal{J}$, containing 
a $1$-morphism that is not isomorphic to any of the identity $1$-morphisms;
\item the cell $\mathcal{J}$ is the maximal two-sided cell with respect to $\leq_J$;
\item the cell $\mathcal{J}$ is idempotent in the sense that $\mathrm{F}\circ \mathrm{G}\neq 0$
for some (possibly equal) $1$-morphisms $\mathrm{F}$ and $\mathrm{G}$ in $\mathcal{J}$.
\end{itemize}
In particular, a simple $2$-category is $\mathcal{J}$-simple in the sense of \cite[Subsection~6.2]{MM2}.

We note that the above definition excludes the very easy situation when the only indecomposable
$1$-morphisms in $\cC$ are those isomorphic to the identity $1$-morphisms. Such $2$-categories are
easy to construct and study, so we will ignore them.

Similarly to \cite[Subsection~6.2]{MM2}, one shows that a finitary $2$-category which satisfies
the last three of the above conditions has a unique simple quotient.

\subsection{Simple fiat strongly regular $2$-categories}\label{s4.2}

As in \cite[Subsection~6.2]{MM1}, a finitary $2$-category $\cC$ is called {\em fiat} provided 
that it has a weak anti-involution $\star$ (reversing both $1$- and $2$-morphisms) and
adjunction morphisms between any pair $(\mathrm{F},\mathrm{F}^{\star})$ of $1$-morphisms. 
A fiat $2$-category $\cC$ is called {\em strongly regular} provided that no two left (right) cells inside the same two-sided cell are comparable, and that the intersection of 
each left and each right cell inside the same two-sided cell consists of exactly one element,
see \cite[Subsection~4.8]{MM1} and \cite[Corollary~19]{KM2}.

$2$-categories which are, at the same time, simple, strongly regular and fiat, were classified
in \cite[Theorem~13]{MM3}. Roughly speaking (up to the structure of the endomorphism algebra 
of identity $1$-morphisms), they are biequivalent to the bicategory of projective 
bimodules for a finite dimensional, weakly symmetric $\Bbbk$-algebra $A$, or, equivalently to the $2$-category $\cC_A$  defined as follows:
Let $A=A_1\times A_2\times\dots\times A_k$ be a decomposition of $A$ into a direct sum of connected 
components. Then
\begin{itemize}
\item objects of $\cC_A$ are $\mathtt{1},\mathtt{2},\dots,\mathtt{k}$, where 
$\mathtt{i}$ should be thought of as a small category $\mathcal{A}_i$ equivalent to $A_i$-mod;
\item $1$-morphisms in $\cC_A(\mathtt{i},\mathtt{j})$ are functors 
from $\mathcal{A}_i$ to $\mathcal{A}_j$ corresponding to tensoring
with bimodules from the additive closure of $A_j\otimes_{\Bbbk}A_i$, with the additional
bimodule $A_i$, if $i=j$;
\item $2$-morphisms are all natural transformations of such functors.
\end{itemize}

It is natural to ask for a description of simple finitary $2$-categories in general. 
Unfortunately, an easy generalisation of \cite[Theorem~13]{MM3} is too much to hope for. 

\subsection{Cell $2$-representations}\label{s4.3}

For any finitary $2$-category $\cC$, we can consider the $2$-category $\cC$-afmod
of {\em finitary $2$-representations of $\cC$}, where
\begin{itemize}
\item objects in $\cC$-afmod are $2$-functors which represent every object in
$\cC$ by a finitary $\Bbbk$-linear category, each $1$-morphism by an additive functor
and each $2$-morphism by a natural transformation of functors;
\item $1$-morphisms in $\cC$-afmod are $2$-natural transformations, see \cite[Subsection~2.3]{MM3};
\item $2$-morphisms in $\cC$-afmod are modifications.
\end{itemize}
An easy example of a $2$-representation is the {\em principal} (Yoneda) $2$-representation
$\mathbf{P}_{\mathtt{i}}:=\cC(\mathtt{i},{}_-)$, defined for each $\mathtt{i}\in\cC$.

If $\mathcal{L}$ is a left cell, then there is $\mathtt{i}\in \cC$ such that all $1$-morphisms 
in $\mathcal{L}$ start at $\mathtt{i}$. In this case we can define the corresponding
{\em cell $2$-representation}  $\mathbf{C}_{\mathcal{L}}$ of $\cC$ as a certain
subquotient of $\mathbf{P}_{\mathtt{i}}$, see \cite[Subsection~6.5]{MM2}, and also
\cite[Section~4]{MM1} for the original definition.

\subsection{Main idea}\label{s4.4}

Let now $\cC$ be a simple finitary $2$-category. Fix a left cell $\mathcal{L}$ inside
the distinguished two-sided cell $\mathcal{J}$ and consider the corresponding cell
$2$-representation $\mathbf{C}_{\mathcal{L}}$. By construction, the fact that 
$\mathcal{J}$ is idempotent implies that $\mathbf{C}_{\mathcal{L}}$ has trivial annihilator.
Therefore $\mathbf{C}_{\mathcal{L}}$ defines a faithful $2$-functor from $\cC$ into some
$2$-category of right exact functors between modules categories of finite dimensional
algebras. The main point of \cite[Theorem~13]{MM3} was to show that, on the level of 
$1$-morphisms in $\mathcal{J}$, this embedding is, in fact, $2$-full and each 
$1$-morphism in $\mathcal{J}$ is represented by a {\em projective functor}, that is
by a functor isomorphic to tensoring with a projective bimodule.

If $\cC$ is fiat (but not necessarily strongly regular), then \cite[Theorem~2]{KMMZ}
still asserts that, in the setup of the previous paragraph, each $1$-morphism in 
$\mathcal{J}$ is  represented by a {\em projective functor}. However, outside
strongly regular situation the $2$-fullness statement is no longer true. An easy
example is given by the small quotient of the $2$-category of Soergel bimodules in
Weyl type $B_2$, see \cite[Subsection~3.2]{KMMZ}. This is a simple, fiat, but not
strongly regular $2$-category. Under $\mathbf{C}_{\mathcal{L}}$, the indecomposable
$1$-morphism corresponding to a simple reflection belonging to  $\mathcal{L}$ acts
non-trivially on two simple modules of the abelianised cell $2$-representation. Therefore
this indecomposable $1$-morphism is represented by a decomposable projective functor,
which means that the representation $2$-functor cannot be $2$-full.

If $\cC$ is not fiat, then $1$-morphisms in $\mathcal{J}$ no longer act as  
projective functors in general. Later on we will use our Theorem~\ref{thm1}
to describe what kind of functors appear in this case.

\subsection{Simple transitive $2$-representations}\label{s4.5}

Cell $2$-representations of finitary $2$-categories have the following two properties:
\begin{itemize}
\item They are {\em transitive} in the sense that, for any pair $X,Y$ of indecomposable
objects in the underlying categories of the representation, there is a $1$-morphism
$\mathrm{F}$ such that $X$ appears (up to isomorphism) as a direct summand in
$\mathrm{F}\,Y$.
\item They are {\em simple} in the sense that the underlying categories do not have any proper
ideals invariant under the action of $\cC$.
\end{itemize}
In general, simple transitive $2$-representations are natural $2$-analogues of simple
modules over $\Bbbk$-algebras, see \cite{MM5,MM6}.

If $\mathbf{M}$ is a simple transitive $2$-representation of $\cC$, then, according to
\cite[Subsection~3.2]{CM}, there is a unique two-sided cell $\mathcal{J}_{\mathbf{M}}$
which is a maximal element (with respect to the two-sided order) in the set of all 
two-sided cells whose $1$-morphisms are not annihilated by $\mathbf{M}$. 
The cell $\mathcal{J}_{\mathbf{M}}$ is called the {\em apex of
$\mathbf{M}$} and is idempotent in the sense that it contains three (possibly equal)
$1$-morphisms $\mathrm{F}$, $\mathrm{G}$ and $\mathrm{H}$ such that 
$\mathrm{H}$ appears, as a summand and up to isomorphism, in $\mathrm{F}\circ \mathrm{G}$.

\subsection{Abelianisation}\label{s4.6}

Instead of {\em additive} $2$-representations of $\cC$ one can also study {\em abelian}
$2$-representations where 
\begin{itemize}
\item each object in $\cC$ is represented by a category equivalent to the module
category for some finite dimensional associative $\Bbbk$-algebra;
\item each $1$-morphism is represented by a right exact functor;
\item each $2$-morphism is represented by a natural transformation of functors.
\end{itemize}
The $2$-category of abelian $2$-representations is denoted $\cC$-mod, see e.g. \cite[Section~4]{MM2}.
There is a natural $2$-functor $\overline{\,\,\cdot\,\,}:\cC\text{-}\mathrm{afmod}\to
\cC\text{-}\mathrm{mod}, \mathbf{M}\mapsto \overline{\mathbf{M}}$, called {\em abelianisation $2$-functor}, see \cite[Subsection~4.2]{MM2}.

\subsection{Main results of this section}\label{s4.7}

\begin{proposition}\label{prop31}
Let $\cC$ be a finitary $2$-category and $\mathbf{M}$ 
a simple transitive $2$-rep\-re\-sen\-ta\-tion of $\cC$. Then, for any $1$-morphism $\mathrm{F}$ in
$\mathcal{J}_{\mathbf{M}}$, the functor $\overline{\mathbf{M}}(\mathrm{F})$  sends any
object to a projective object.
\end{proposition}

\proof
The claim follows directly from the proof of the first part of \cite[Theorem~2]{KMMZ}
as this proof does not use fiatness of $\cC$.
\endproof

Let $A$ be a finite dimensional associative $\Bbbk$-algebra with a fixed decomposition 
$A=A_1\times A_2\times\dots\times A_k$ into a product of (not necessarily connected) components.
For each $i=1,2,\dots,k$ fix a small category $\mathcal{C}_i$ equivalent to 
$A_i$-mod and a right $A_i$-module $N_i$. Let $\mathcal{C}=\{\mathcal{C}_i\}$ and
$N=\{N_i\}$. Define the weakly finitary $2$-category $\cC_{A,\mathcal{C},N}$ as follows:
\begin{itemize}
\item  The objects of $\cC_{A,\mathcal{C},N}$ are $\mathtt{1},\mathtt{2},\dots,\mathtt{k}$,
where $\mathtt{i}$ should be thought of as $\mathcal{C}_i$;
\item $1$-morphisms in $\cC_{A,\mathcal{C},N}(\mathtt{i},\mathtt{j})$ are all functors
from $\mathcal{C}_i$ to $\mathcal{C}_j$ which are isomorphic to tensoring with 
$A_j$-$A_i$-bimodules in the additive closure of $A_j\otimes_{\Bbbk}N_i$ and, additionally,
the $A_i$-$A_i$-bimodule $A_i$, if $i=j$;
\item $2$-morphisms are natural transformations of functors.
\end{itemize}

The main result of this section is the following statement which, roughly speaking, says that
all simple finitary $2$-categories are $2$-subcategories of the categories of the
form $\cC_{A,\mathcal{C},N}$.

\begin{theorem}\label{thm33}
Let $\cC$ be a simple finitary $2$-category. Then there are $A$, $\mathcal{C}$ and $N$ as above
and a faithful $2$-functor from $\cC$ to $\cC_{A,\mathcal{C},N}$.
\end{theorem}

\proof
Let $\cC$ be a simple finitary $2$-category. For a left cell $\mathcal{L}$ in $\mathcal{J}$,
consider the left cell $2$-representation $\mathbf{C}_{\mathcal{L}}$ of $\cC$.
By Proposition~\ref{prop31}, for any $1$-morphism $\mathrm{F}$ in $\mathcal{J}$
from $\mathtt{i}$ to $\mathtt{j}$,
the functor $\overline{\mathbf{C}_{\mathcal{L}}}(\mathrm{F})$ maps any object
of $\overline{\mathbf{C}_{\mathcal{L}}}(\mathtt{i})$ to a projective object
of $\overline{\mathbf{C}_{\mathcal{L}}}(\mathtt{j})$. 

For each $\mathtt{i}$, let $A_{\mathtt{i}}$ denote the underlying algebra of 
$\overline{\mathbf{C}_{\mathcal{L}}}(\mathtt{i})$. We note that 
$A_{\mathtt{i}}$ does not have to be connected. By Theorem~\ref{thm1},
there is a right $A_{\mathtt{i}}$-module $N_i$ such that any $1$-morphism in 
$\mathcal{J}$ from $\mathtt{i}$ to any $\mathtt{j}$ is represented, via 
$\overline{\mathbf{C}_{\mathcal{L}}}$, by a functor isomorphic to tensoring 
with a bimodule of the form $A_j\otimes_{\Bbbk}N_i$ (and, additionally, $A_i$,
if $i=j$). 

Let $\displaystyle A:=\prod_{\mathtt{i}}A_{\mathtt{i}}$. 
Since $\mathcal{J}$ is idempotent, it coincides with the apex of $\mathbf{C}_{\mathcal{L}}$.
Hence, thanks to simplicity of $\cC$, the representation $2$-functor 
$\mathbf{C}_{\mathcal{L}}$ is faithful on the level of $2$-morphisms
and thus provides an embedding of $\cC$ into the weakly finitary category 
$\cC_{A,\{\overline{\mathbf{C}_{\mathcal{L}}}(\mathtt{i})\},\{N_{\mathtt{i}}\}}$. 
This completes the proof.
\endproof

We note that, usually, the embedding of $\cC$ given by Theorem~\ref{thm33} will not be $2$-full.
Furthermore, $A$, $\mathcal{C}$ and $N$ in the formulation of Theorem~\ref{thm33}
are not uniquely defined, even up to isomorphism/equivalence.

\section{Decategorification}\label{snew}

Let $\cC:=\cC_{A,\mathcal{C},N}$ be as in Section~\ref{s4.7}. Let $P_1,\dots,P_r$ be a complete list of projective indecomposable modules for $A$, and $N_1,\dots, N_s$ a complete list of elements in $N$. Then a complete list of indecomposable $1$-morphisms in $\cC_{A,\mathcal{C},N}$ is given by $\mathbbm{1}_\mathtt{i}$, $\mathtt{i}\in\{\mathtt{1}, \dots, \mathtt{k}\}$, and $\{\mathrm{F}_{i,j}\}$ where  $\mathrm{F}_{i,j}$ is a functor isomorphic to tensoring with $P_i\otimes_\Bbbk N_j$. 
The structure of a multisemigroup with finite multiplicities (in the sense of \cite{Fo}) 
on $\mathcal{S}[\cC]\cap \mathcal{J}$ is given by $[\mathrm{F}_{i,j}][\mathrm{F}_{l,t}] = \dim (N_je_l) [\mathrm{F}_{i,t}]$.

\begin{proposition}\label{decat}
Suppose $\mathrm{add}(N_1\oplus N_2\oplus \dots \oplus N_s)\cong \mathrm{add}(A)$. 
\begin{enumerate}[$($i$)$]
\item\label{filt} The algebra $\Bbbk\otimes_{\mathbb{Z}} \mathbb{Z}\mathcal{S}[\cC]$ has a filtration by two-sided ideals such that the lowest ideal is given by the span of $\{[\mathrm{F}_{i,j}]| i,j=1,\dots,r\}$, and the remaining subquotients are spanned by $[\mathbbm{1}_\mathtt{i}]$, taken in any order.

\item\label{swich} The ideal $J$ spanned by $\{[\mathrm{F}_{i,j}]| i,j=1,\dots,r\}$ in $\Bbbk\otimes_{\mathbb{Z}} \mathbb{Z}\mathcal{S}[\cC]$ is, in the terminology of \cite[Definition 3.3]{KX}, 
isomorphic to the swich algebra of the algebra of $r\times r$-matrices with respect to the matrix $(\dim(e_jAe_l))_{j,l=1}^r$.

\item\label{cell} If $A$ is weakly symmetric, the filtration in \eqref{filt} is a cell filtration, where the involution is given by the action of ${}^*$ on $\mathcal{S}[\cC]$, which corresponds to interchanging the two subscripts on $\mathrm{F}_{i,j}$.
\end{enumerate}
\end{proposition}

\begin{proof}
\eqref{filt} and \eqref{swich} follow directly from the definitions.

To prove \eqref{cell}, assume $A$ is weakly symmetric, so $\cC$ is fiat. The involution ${}^*$ induces an involution on $\Bbbk\otimes_{\mathbb{Z}} \mathbb{Z}\mathcal{S}[\cC]$, fixing the $[\mathbbm{1}_\mathtt{i}]$, and sending $ [\mathrm{F}_{i,j}]$ to $[\mathrm{F}_{j,i}]$ and hence, under the isomorphism to the swich algebra from \eqref{swich}, acts as matrix transposition. Since $\dim e_iAe_j = \dim e_jAe_i$, it is easy to check that, for $V$ the $\Bbbk$-span of $[\mathrm{F}_{i,1}]$, where $i=1,2,\dots,r$, the morphism $\alpha: \, J \to V\otimes V, \quad \mathrm{F}_{j,i} \mapsto v_i\otimes v_j$, defines the structure of a cell ideal on $J$ (cf. \cite[Definition 3.2]{KX1}).

Quotienting out by $J$, all remaining subquotients in the ideal filtration are one dimensional and idempotent, and hence cell ideals.
\end{proof}

\begin{corollary}
Let $\cC$ be a fiat $2$-category such that all two-sided cells are strongly regular. Then $\Bbbk\otimes_{\mathbb{Z}} \mathbb{Z}\mathcal{S}[\cC]$ is a quasi-hereditary algebra with simple-preserving duality.
\end{corollary}

\begin{proof}
By induction with respect to the two-sided order, it follows immediately from 
Proposition \ref{decat}\eqref{cell} that $\Bbbk\otimes_{\mathbb{Z}} \mathbb{Z}\mathcal{S}[\cC]$ 
is a cellular algebra. Since each cell contains an idempotent, it is indeed quasi-hereditary.
\end{proof}

\section{Application: simple transitive $2$-representations of projective bimodules
for $\Bbbk[x,y]/(x^2,y^2,xy)$}\label{s5}

\subsection{Classification}\label{s5.1}

In this section we consider the problem of classification of simple transitive 
$2$-representations of the simple $2$-category $\cC_A$ of projective bimodules for the 
$3$-dimensional algebra $A=\Bbbk[x,y]/(x^2,y^2,xy)$. As $A$ is connected,
the $2$-category $\cC_A$ has only one object. We call this object $\mathtt{i}$.
The main result of this section is the following.

\begin{theorem}\label{thm51}
Every simple transitive  $2$-representations of $\cC_A$, for the algebra
$A=\Bbbk[x,y]/(x^2,y^2,xy)$,
is equivalent to a cell $2$-representation.
\end{theorem}

Under the additional assumption that all $1$-morphisms in $\cC_A$ are represented by 
exact functors, the claim of Theorem~\ref{thm51} is proved in \cite[Proposition~19]{MM5}.
Therefore, to prove Theorem~\ref{thm51}, we just have to show that, in every simple 
transitive  $2$-representations of $\cC_A$, all $1$-morphisms in $\cC_A$ are indeed 
represented by exact functors. The rest of this section is devoted to the proof of
Theorem~\ref{thm51}.

\subsection{Combinatorics of the action}\label{s5.2}

We fix a simple transitive $2$-representation $\mathbf{M}$ of $\cC_A$. Let
$X_1,X_2,\dots,X_k$ be a list of representatives of isomorphism classes of indecomposable
objects in $\mathbf{M}(\mathtt{i})$. Let $\mathrm{F}$ be an indecomposable $1$-morphism
in $\cC_A$ which is isomorphic to tensoring with $A\otimes_{\Bbbk}A$. Note that 
\begin{equation}\label{eq54}
\mathrm{F}\circ \mathrm{F}\cong \mathrm{F}^{\oplus 3}.
\end{equation}
For $i,j=1,2,\dots,k$, let 
$m_{i,j}$ denote the multiplicity of $X_i$ in $\mathrm{F}\, X_j$. Then 
$M=(m_{i,j})_{i,j=1}^k\in\mathrm{Mat}_{k\times k}(\mathbb{Z}_{\geq0})$, moreover,
all $m_{i,j}$ are positive due to transitivity of $\mathbf{M}$, see \cite[Subsection~7.1]{MM5}
for details.

From \eqref{eq54}, it follows that, $M^2=3M$.
Therefore, up to permutation of the $X_i$'s, $M$ is one of the following matrices 
(again, see \cite[Subsection~7.1]{MM5} and \cite[Proposition~10]{MZ} for details):
\begin{displaymath}
M_1=\left(3\right),\quad
M_2=\left(\begin{array}{cc}2&2\\1&1\end{array}\right),\quad
M_3=\left(\begin{array}{cc}2&1\\2&1\end{array}\right),\quad
M_4=\left(\begin{array}{ccc}1&1&1\\1&1&1\\1&1&1\end{array}\right).
\end{displaymath}

\subsection{General part of the proof of Theorem~\ref{thm51}}\label{s5.3}

Let $B$ be an associative algebra such that $\overline{\mathbf{M}}(\mathtt{i})\cong B$-mod, then we have
${\mathbf{M}}(\mathtt{i})\cong B$-proj. Let $e_1,e_2,\dots,e_k$
be primitive idempotents of $B$ corresponding to the objects $X_1,X_2,\dots,X_k$ above.
Note that $k\leq 3$ by Subsection~\ref{s5.2}. For $i=1,2,\dots,k$, we denote by $L_i$
the simple top of $X_i$ in $\overline{\mathbf{M}}(\mathtt{i})$. We also denote by 
$L'_i$ the corresponding  {\em right} simple $B$-module, for $i=1,2,\dots,k$.

By Proposition~\ref{prop31}, the functor
$\overline{\mathbf{M}}(\mathrm{F})$ sends any object in $\overline{\mathbf{M}}(\mathtt{i})$
to a projective object in $\overline{\mathbf{M}}(\mathtt{i})$. Hence, by Theorem~\ref{thm1},
there are right $B$-modules $N_1,N_2\dots,N_k$ such that $\overline{\mathbf{M}}(\mathrm{F})$ 
is isomorphic to tensoring with the $B$-$B$-bimodule
\begin{displaymath}
Be_1\otimes_{\Bbbk}N_1\oplus  
Be_2\otimes_{\Bbbk}N_2\oplus\dots\oplus  
Be_k\otimes_{\Bbbk}N_k.
\end{displaymath}
Note that, for $i,j=1,2,\dots,k$, we have 
\begin{displaymath}
Be_i\otimes_{\Bbbk}N_i\otimes_B Be_j\cong Be_i^{\oplus \dim_{\Bbbk}(N_i\otimes_BBe_j)}
\end{displaymath}
and hence 
\begin{displaymath}
m_{i,j}=\dim_{\Bbbk}(N_i\otimes_BBe_j)=\dim_{\Bbbk}(N_ie_j)=[N_i:L'_j].
\end{displaymath}

Next we observe that the right $B$-module $N=N_1\oplus N_2\oplus\dots\oplus N_k$ is faithful.
Indeed, if this were not the case, the annihilator of this module would be, due to the form of
the functor $\overline{\mathbf{M}}(\mathrm{F})$, annihilated by $\overline{\mathbf{M}}(\mathrm{F})$ 
and hence would generate a non-zero proper $\cC_A$-stable ideal of $\mathbf{M}$. 
This is, however, not possible thanks to  simplicity of $\mathbf{M}$. Consequently, as the sum 
of entries in each row of the above matrices is at most $4$, the Loewy length of $N$ is at most $4$, and
we have that $\mathrm{Rad}(B)^4=0$.

As mentioned above, we just need to show that $\overline{\mathbf{M}}(\mathrm{F})$ is exact or,
equivalently, that $N$ is projective. If $B$ is semi-simple, then $N$ is automatically projective,
so in what follows we may assume that $\mathrm{Rad}(B)\neq 0$.

Finally, we will need the following.

\begin{lemma}\label{lem541}
The $2$-functor $\overline{\mathbf{M}}$ induces an embedding of the algebra $A$ into the center 
$Z(B)$ of the algebra $B$. 
\end{lemma}

\begin{proof}
The $2$-functor $\overline{\mathbf{M}}$ gives a map from 
$A\cong\mathrm{End}_{\ccC_A}(\mathbbm{1})$ to the endomorphism algebra of
the $B$-$B$-bimodule $B$. The latter is isomorphic to $Z(B)$ and injectivity of
this map follows from simplicity of $\cC_A$.
\end{proof}

Now we have to go into a case-by-case analysis.

\subsection{Case~1: $M=M_1$}\label{s5.4}

In this case we have $k=1$ and $N=N_1$ has dimension $3$. If $\mathrm{Rad}(B)^2\neq 0$, then
$N$ must be uniserial and hence $B\cong\Bbbk[x]/(x^3)$. This means that $N$ is projective and we are done.

If $\mathrm{Rad}(B)^2=0$, then we have two possibilities. The first one is $B\cong\Bbbk[x]/(x^2)$
and $N=B\oplus \Bbbk$, which immediately contradicts Lemma~\ref{lem541}.
The second possibility which is left is $B\cong A$. In this case, as $N$ is faithful, it must be either the 
projective indecomposable or the injective indecomposable $B$-module. In the projective case we are done, 
so let us assume that $N$ is injective. To get a contradiction in this case we will need to do more subtle computations. 

We have  $A=B=\Bbbk[x,y]/(x^2,y^2,xy)$. Then the $A$-$A$-bimodule $A\otimes_{\Bbbk} A$ has the following basis:
\begin{displaymath}
1\otimes 1,\,1\otimes x,\,1\otimes y,\, x\otimes 1,\, y\otimes 1,\,x\otimes x,\,x\otimes y,\, y\otimes x,\, 
y\otimes y. 
\end{displaymath}
The $A$-$A$-bimodule $A\otimes_{\Bbbk} N\cong A\otimes_{\Bbbk} A^{*}$ has the following basis:
\begin{displaymath}
1\otimes x^*,\,1\otimes y^*,\,1\otimes 1^*,\, x\otimes x^*,\, x\otimes y^*,\,y\otimes x^*,\,y\otimes y^*,\, 
x\otimes 1^*, y\otimes 1^*. 
\end{displaymath}
Now it is easy to check that $\dim_{\Bbbk}\mathrm{Hom}_{A\text{-}A}(A,A\otimes_{\Bbbk}A)=4$, where 
the generator $1$ of $A$ can be mapped to any of the elements
\begin{displaymath}
x\otimes x,\,x\otimes y,\, y\otimes x,\, y\otimes y.
\end{displaymath}
At the same time, $\dim_{\Bbbk}\mathrm{Hom}_{A\text{-}A}(A,A\otimes_{\Bbbk}A^*)=3$, where 
the generator $1$ of $A$ can be mapped to any of the elements
\begin{displaymath}
1\otimes 1^*+x\otimes x^*+y\otimes y^*,\, x\otimes 1^*,\, y\otimes 1^*.
\end{displaymath} 
As $\cC_A$ is simple, $\mathrm{Hom}_{A\text{-}A}(A,A\otimes_{\Bbbk}A)$ should be embeddable into
the homomorphism spaces between the functors $\overline{\mathbf{M}}(\mathbbm{1}_{\mathtt{i}})$ and 
$\overline{\mathbf{M}}(\mathrm{F})$, but the above calculation shows that this is not possible.
This completes the proof in Case~1.

\subsection{Case~2: $M=M_2$}\label{s5.5}

In this case $k=2$, $[N_1:L'_i]=2$ and $[N_2:L'_i]=1$, for $i=1,2$. 
The endomorphism algebra of the multiplicity-free module $N_2$ is a direct sum of copies of $\Bbbk$.
As $N$ is faithful, it follows that $A$ embeds into $\mathrm{End}_B(N_1)$. Both simples
appear in $N_1$ with multiplicity $2$. This implies that the image of $Z(B)$ in the endomorphism
algebra of each of the two indecomposable projective $B$-modules has dimension at most $2$. 
Therefore both these images must have dimension $2$ and the corresponding kernels must be different.

As $A$ embeds into the endomorphism algebra of $N_1$, we have
$\dim_{\Bbbk}(\mathrm{End}_B(N_1))\geq 3$. In particular, $N_1$ cannot have simple top or socle for in that case
the dimension of its endomorphism algebra would be at most $2$, that is the multiplicity of the top
(or socle) simple module in $N_1$. As $N_1$ cannot be semi-simple (since $A$ cannot be embedded into a
direct sum of copies of $\Bbbk$), it follows that $N_1$ has Loewy length two with both top
and socle isomorphic to $L'_1\oplus L'_2$. Consequently, $\mathrm{Rad}(B)^2=0$.

It follows that the quiver of $B$ has one loop at each vertex. 
Furthermore, we could also have at most two arrows in each direction
between the two different vertices (and, in total, at most three such arrows). 
Thus, if $B$ is not connected, then 
$B=Z(B)\cong \Bbbk[x]/(x^2)\oplus \Bbbk[y]/(y^2)$. If $B$ is connected, then $Z(B)\cong A$.

Consider first the case when  $B$ is not connected. In this case the above discussion implies
$N_1\cong B$ while $N_2\cong \Bbbk_{x}\oplus \Bbbk_{y}$, the direct sum of two non-isomorphic 
simple $B$-modules corresponding to the $x$- and  $y$-components of $B$, respectively. Setting
$D_x:=\Bbbk[x]/(x^2)$ and $D_y:=\Bbbk[y]/(y^2)$, we obtain that $\overline{\mathbf{M}}(\mathrm{F})$
is represented by the bimodule
\begin{equation}\label{eq58}
(D_x\otimes_{\Bbbk}D_x)\oplus (D_x\otimes_{\Bbbk}D_y)\oplus (D_y\otimes_{\Bbbk}\Bbbk_x)\oplus 
(D_y\otimes_{\Bbbk}\Bbbk_y). 
\end{equation}
The dimension of the homomorphism space from the $B$-$B$-bimodule $B$ to the bimodule 
in \eqref{eq58} is $3$, where two
linearly independent homomorphisms come from homomorphisms from $D_x$ to  $D_x\otimes_{\Bbbk}D_x$
and one more linearly independent homomorphism comes from the map from $D_y$ to $D_y\otimes_{\Bbbk}\Bbbk_y$.
As mentioned above, $\dim_{\Bbbk}\mathrm{Hom}_{A\text{-}A}(A,A\otimes_{\Bbbk}A)=4$ and
we get a contradiction.

If $B$ is connected, the analogue of the bimodule \eqref{eq58} will look more complicated.
However, the new quiver will be obtained from the case of disconnected $B$ by adding some
new arrows between different vertices. The restriction of the center of this new $B$ to each
vertex still coincides with the case considered in the previous paragraph. Therefore any
homomorphism from $B$ to this new bimodule will restrict, by forgetting all new arrows,
to a homomorphism from the previous paragraph. As generators for our bimodules remain the
same, this restriction map is injective. Of course, new arrows might come with new conditions,
so there is no obvious guarantee that the restriction map is surjective. In any case, the 
dimension of the homomorphism space from our new $B$ to this new bimodule will be at most $3$,
which again leads to a contradiction. This completes the proof in Case~2.

\subsection{Case~3: $M=M_3$}\label{s5.6}

In this case $k=2$, $[N_i:L'_1]=2$ and $[N_i:L'_2]=1$, for $i=1,2$. 
As $N$ is faithful, we obviously have $\mathrm{Rad}(B)^3=0$. Moreover, as $N$ is faithful and 
$L'_2$ appears with multiplicity $1$ in both $N_1$ and $N_2$, it follows that 
$\mathrm{End}(P_2)\cong\Bbbk$ and hence Lemma~\ref{lem541} gives an embedding
of $A$ into $\mathrm{End}(P_1)$. As $N$ has only $2$ direct summands and 
$L_1$ appears with multiplicity $2$ in both, it follows that $A\cong \mathrm{End}(P_1)$ and
that the generators $x$ and $y$ of $A$ can be chosen such that 
$\mathrm{End}(N_1e_1)\cong A/(x)$ under the natural map from $A$ to $\mathrm{End}(N_1e_1)$.

The $2$-functor $\overline{\mathbf{M}}$ induces a map from the $4$-dimensional 
space $\mathrm{Hom}_{A\text{-}A}(A,A\otimes_{\Bbbk}A)$ to the space of $B$-$B$-bimodule
homomorphisms from $B$ to $P_1\otimes_{\Bbbk}N_1\oplus P_2\otimes_{\Bbbk}N_2$. Due 
to the previous paragraph, the image of the generator $1\in B$ belongs to
\begin{displaymath}
e_1P_1\otimes_{\Bbbk}N_1e_1\oplus e_2P_2\otimes_{\Bbbk}N_2e_2.
\end{displaymath}
The space $e_2P_2\otimes_{\Bbbk}N_2e_2$ has dimension $1$ by the above. The space 
$e_1P_1\otimes_{\Bbbk}N_1e_1$ can be identified with $A\otimes_{\Bbbk}\left(A/(x)\right)$.
Even on the level of $A$-$A$-bimodules, the image of $1$ in  $A\otimes_{\Bbbk}A/(x)$
has to belong to the $2$-dimensional subspace
generated by $x\otimes y$ and $y\otimes y$. Therefore we obtain an injection of 
a $4$-dimensional space into a space of dimension at most $3$, a contradiction. 
The proof of Case~3 is now complete.

\subsection{Case~4: $M=M_4$}\label{s5.7}

In this case $k=3$ and $[N_i:L'_j]=1$, for all $i,j=1,2,3$. In particular,
all  $N_i$'s are multiplicity free and hence the endomorphism algebra of
each $N_i$ is a direct sum of copies of $\Bbbk$.
Thus $\overline{\mathbf{M}}$ cannot embed $A$ into $Z(B)$, contradicting 
Lemma~\ref{lem541}. This completes the proof in Case~4 and thus of the whole Theorem~\ref{thm51}.

\vspace{5mm}

\noindent
Vo.~Ma.: Department of Mathematics, Uppsala University, Box. 480,
SE-75106, Uppsala, SWEDEN, email: {\tt mazor\symbol{64}math.uu.se}

\noindent
Va.~Mi.: School of Mathematics, University of East Anglia,
Norwich NR4 7TJ, UK,  {\tt v.miemietz\symbol{64}uea.ac.uk}

\noindent
X.~Z.: Department of Mathematics, Uppsala University, Box. 480,
SE-75106, Uppsala, SWEDEN, email: {\tt xiaoting.zhang\symbol{64}math.uu.se}

\end{document}